\documentclass[11pt]{amsart}

\usepackage{graphicx}

\usepackage{amssymb}

\usepackage{latexsym}

\newtheorem{thm}{Theorem}[section]

\newtheorem{pro}[thm]{Proposition}

\author[O. Goubet]{O. Goubet}

\title[ Zhidkov spaces]{Two remarks on solutions of Gross-Pitaevskii equations on Zhidkov spaces}

\address[Olivier Goubet]{LAMFA CNRS UMR 6140 \\
Universit\'e de Picardie Jules Verne \\
33 rue Saint-Leu 80039 Amiens cedex.}

\email{olivier.goubet@u-picardie.fr}

\begin{document}

\begin{abstract}
 We consider the so-called Gross-Pitaevskii equations supplemented
 with non-standard boundary conditions. We prove two mathematical
 results concerned with the initial value problem for these
 equations in Zhidkov spaces.
\end{abstract}


\subjclass[2000]{}

\date{6 july 2005}

\maketitle

\section{Introduction}

\subsection{Setting the problem}

This short article is concerned with the so-called Gross-Pitaevskii
equations
\begin{equation}\label{gp}
iu_t+\Delta u + u(1-|u|^2)=0,
\end{equation}
\noindent supplemented with non-standard boundary conditions that
read $|u(t,x)|\rightarrow 1$ as $||x||\rightarrow +\infty$. We also
supplement this equation with initial condition $u_0$ that will be
specified in the sequel. Here the unknown $u$ maps
$\mathbb{R}_t\times \mathbb{R}_x^D$ into $\mathbb{C}$. These
equations with this non-standard boundary conditions occur in
several physical contexts, as the Bose-Einstein condensation for
suprafluids (see \cite{Gallo1} and the references therein). The
mathematical study of solitary waves for these equations was
initiated in the pioneering work \cite{BethuelSaut}. Here we are
interested in two questions related to the Cauchy problem for
(\ref{gp}). Throughout this article, we follow the framework
developed by C. Gallo (see \cite{Gallo1}), that extends to the
multi-dimensional case the work of P. Zhidkov (see \cite{Zhidkov}
and the references therein).

We now introduce the Zhidkov spaces as, for $k\geq 1$ integer

\begin{equation}\label{Z}
X^k(\mathbb{R}^D)=\{ u \in L^\infty(\mathbb{R}^D)\cap
UC(\mathbb{R}^D); \nabla u \in H^{k-1}(\mathbb{R}^D)\}.
\end{equation}

\noindent Here $UC(\mathbb{R}^D)$ denotes the space of uniformly
continuous functions and $H^{k-1}(\mathbb{R}^D)$ is the usual
Sobolev space. This choice is suitable for the above boundary
conditions. The norm on $X^k(\mathbb{R}^D)$ is

\begin{equation}\label{znorm}
||u||_{X^k(\mathbb{R}^D)}=||u||_{L^\infty(\mathbb{R}^D)}+||\nabla
u||_{H^{k-1}(\mathbb{R}^D)}.
\end{equation}

\noindent C. Gallo has proved the following result
\begin{thm}
Assume $k> \frac{D}{2}$ and $u_0\in X^k(\mathbb{R}^D)$. Then the
initial value problem for Gross-Pitaevskii equation is locally well
posed in $X^k(\mathbb{R}^D)$. If moreover $u_0\in
X^{k+2}(\mathbb{R}^D)$, then $u(t)\in X^{k+2}(\mathbb{R}^D)$.
\end{thm}

\noindent On the other word, due to the particular form of the
boundary conditions, we define the energy associated to our problem
as

\begin{equation}\label{energy}
E(u)=\frac{1}{2}\int_{\mathbb{R}^D}|\nabla
u|^2dx+\frac{1}{4}\int_{\mathbb{R}^D} (1-|u|^2)^2dx.
\end{equation}

\noindent C. Gallo also proved that for $E(u_0)< +\infty$ and for
$D\leq 2$, then the energy $E(u(t))$ is conserved along the
trajectories. As a consequence of this result, he has established
that the one-dimensional initial value problem is globally
well-posed in $X^1(\mathbb{R})$.

\noindent In this short article, our first result shows the
persistence of the energy is valid on {\it any} dimension, i.e

\begin{pro}\label{persistence}
For any $D\geq 2$, if $u_0$ in $X^k(\mathbb{R}^D)$ ($k>\frac{D}{2}$) is
such that $E(u_0)< +\infty$, then the energy $E(u(t))$ is conserved
along the trajectories.
\end{pro}

\noindent The proof relies on the study of the growth rate of the
energy along annuli in $\mathbb{R}^D$. 

Our second result states as follows

\begin{pro}\label{global}
The initial value problem for Gross-Pitaevskii equations is globally
well-posed in $X^2(\mathbb{R}^2)$ for initial data $u_0$ that is in
this space and that has finite energy, i.e $E(u_0)<+\infty$.
\end{pro}

\noindent The proof relies on some Brezis-Gallouet inequality
similar to those that appeared in \cite{BrezisGallouet}.

\noindent After this work was completed, we learned that P. G\'erard
(see \cite{Gerard}) proved the global well-posedness for these
equations, in the cases $D=2$ or $D=3$, in the energy space

\begin{equation}\label{energyspace}
E= \{ u \in H^1_{loc}(\mathbb{R}^D); E(u) < +\infty \}.
\end{equation}

\noindent See also \cite{Gallo2} for where the author uses similar
framework to handle more general nonlinearities and problem in
exterior domains. For the sake of completeness, we would like to
point out also the articles \cite{Ginibre1}, \cite{Ginibre2} where
the authors solve a parabolic regularization of the Gross-Pitaevskii
equations, the so-called complex Ginzburg-Landau equations, in local Sobolev
spaces.

\noindent The next section is devoted to proving Proposition
\ref{persistence}. We complete this article by the proof of
Proposition \ref{global} in a third section.

\section{Persistence of the energy}

\subsection{Splitting of the energy}

Consider a function $u$ that has finite energy $E(u)< +\infty$. Set
\begin{equation}\label{annulus}
C_j=\{ x\in \mathbb{R}^D, j\leq |x| <j+1\}.
\end{equation}
\noindent The kinetic energy of $u$ expands as follows
\begin{equation}\label{kineticexpansion}
\begin{split}
\frac{1}{2}\int_{\mathbb{R}^D}|\nabla u|^2dx= \sum_{j=0}^{+\infty}
k_j,\\ k_j=\frac{1}{2}\int_{C_j}|\nabla u|^2dx.
\end{split}
\end{equation}
\noindent The potential energy reads also
\begin{equation}\label{potentialexpansion}
\begin{split}
\frac{1}{4}\int_{\mathbb{R}^D} (1-|u|^2)^2dx= \sum_{j=0}^{+\infty} p_j,\\
p_j=\frac{1}{4}\int_{C_j}(1-|u|^2)^2dx.
\end{split}
\end{equation}

\noindent Consider now a function $u(t,x)$ that is solution to
Gross-Pitaevskii equations and that starts from $u_0$ that has
finite energy. Consider any $T$ such that $u(t)$ belongs to Zhidkov
space $X^2(\mathbb{R}^D)$ for $|t|\leq T$. Obviously the kinetic
part of the energy of $u(t)$ is bounded by the Zhidkov norm. We will
prove below that the potential energy is also bounded, arguing by
contradiction on the growth rate of $j\rightarrow p_j$.

\subsection{Growth rate of the potential energy}

In this section $C$ will be a constant that depends on
$\sup_{|t|\leq T} ||u(t)||_{X^2(\mathbb{R}^D)}$ and that may vary
from one line to one another. On the one hand, we have
\begin{equation}\label{1}
p_j \leq \frac{1}{4}(1+||u||_{L^\infty}^2)^2 {\rm vol}(C_j) \leq C
j^{D-1},
\end{equation}
\noindent where ${\rm vol}(C_j) $ denotes the $D$ dimensional volume
of the annulus.

On the other hand, we introduce $\theta_j$ a radially symmetric
function on $\mathbb{R}^D$ such that $\theta_j =1$ if $|x|\leq j$,
$\theta_j=0$ if $|x|\geq j+1$, and $|\nabla \theta_j|\leq 2$
everywhere. Multiply (\ref{gp}) by $\bar{u}_t\theta_j$ and integrate
the real part of the resulting equation over $\mathbb{R}^D$ to
obtain
\begin{equation}\label{2}
\frac{d}{dt}\left( \int_{\mathbb{R}^D} (\frac{1}{2}|\nabla
u|^2+\frac{1}{4}(1-|u|^2)^2)\theta_jdx \right)=-{\rm Re}\left(
\int_{\mathbb{R}^D} \bar{u}_t \nabla u.\nabla \theta_j dx \right).
\end{equation}

\noindent We majorize the r.h.s of this equality as follows

\begin{equation}\label{3}
|{\rm Re}\left( \int_{\mathbb{R}^D} \bar{u}_t \nabla u.\nabla
\theta_j dx \right)| \leq 2 ||u_t||_{L^2(C_j)}||\nabla
u||_{L^2(C_j)}.
\end{equation}
\noindent Introducing $K_j=\int_0^T k_j dt$ and $P_j=\int_0^T p_j
dt$, we then infer from (\ref{2})-(\ref{3}) (integrating in time and using
Cauchy-Schwarz inequality)

\begin{equation}\label{4}
\sum_{l=0}^{j-1} (k_l+p_l)\leq E(u_0)+2 \
(\int_0^T||u_t||^2_{L^2(C_j)}dt)^{1/2}(\int_0^T||\nabla
u||^2_{L^2(C_j)}dt)^{1/2}.
\end{equation}

\noindent Going back to the equation (\ref{gp}), we also have

\begin{equation}\label{5}
||u_t||_{L^2(C_j)}\leq ||\Delta
u||_{L^2(C_j)}+||u||_{L^\infty(C_j)}||1-|u|^2||_{L^2(C_j)}\leq
C(1+\sqrt{p_j}).
\end{equation}

\noindent Therefore

\begin{equation}\label{6}
\sum_{l=0}^{j-1} (k_l+p_l)\leq E(u_0)+ C \sqrt{K_j}(1+\sqrt{P_j}).
\end{equation}
\noindent Integrating once more in time, we then obtain
\begin{equation}\label{7}
\sum_{l=0}^{j-1} (K_l+P_l)\leq TE(u_0)+ CT\sqrt{K_j}(1+\sqrt{P_j}).
\end{equation}
\noindent We aim to prove that $\sum_{l=0}^{+\infty}P_l<+\infty$.
Let us argue by contradiction; we now pretend that $\sum_{l=0}^{+\infty}P_l=+\infty$.
Since $K_j \rightarrow 0$, we then infer from (\ref{7}) that
$P_j\rightarrow +\infty$. We now introduce
$a=TE(u_0)+2CT\sup_j(\sqrt{K_j})$ and $Q_j=\max(1,P_j)$. We then infer from
(\ref{7}), dropping some unnecessary terms, that $\frac{Q_{j-1}}{a^2}\leq \sqrt{\frac{Q_j}{a^2}}$.

\noindent There exists $j_0$ such that for $j\geq j_0$
then $Q_j\geq 2a^2$. Therefore
$Q_j\geq a^2 2^{2^{j-j_0}}$. 
This contradicts (\ref{1}). Therefore
$\sum_{l=0}^{+\infty}P_l<+\infty$. Going back to (\ref{6}) and
letting $j\rightarrow +\infty$, we thus obtain
\begin{equation}
E(u(t))\leq E(u_0).
\end{equation}

\noindent Since we can go backward in time in the equation, the
reverse inequality is also valid. This completes the proof of
Proposition \ref{persistence}.

\section{Global existence result in the two-dimensional case}

In this section $c$ is a numerical constant that may vary from one
line to one another.

\subsection{Brezis-Gallou\"et type inequality}

We first state and prove

\begin{pro}\label{bg}
Consider a function $u$ that belongs to $X^2(\mathbb{R}^2)$ and such
that $E(u)< +\infty$. There exists a constant $c$ that is
independent of $u$ such that
$$ ||u||_{L^\infty(\mathbb{R}^2)} \leq c(1+\sqrt{E(u)})\left(1+\log
(1+||\Delta u||^2_{L^2(\mathbb{R}^2)})\right)^{1/2}.$$
\end{pro}

\noindent Proof: Consider first a smooth function $\phi(x)$ in
$H^2(\mathbb{R}^2)$. Consider $R>0$. We have

\begin{equation}\label{14}
\begin{split}
||\hat{\phi}||_{L^1(\mathbb{R}^2_\xi)} \leq c\left(
||\phi||_{H^1(\mathbb{R}^2)} (\int_{|\xi|\leq
R}\frac{d\xi}{1+|\xi|^2})^{1/2} + ||\Delta
\phi||_{L^2(\mathbb{R}^2)}
(\int_{|\xi|\geq R}\frac{d\xi}{|\xi|^4})^{1/2}\right)\leq \\
c(||\phi||_{H^1(\mathbb{R}^2)} \log(1+R^2)^{1/2}+||\Delta
\phi||_{L^2(\mathbb{R}^2)}R^{-1}.
\end{split}
\end{equation}

Consider $u$ in $\cap_{k\geq 1} X^k(\mathbb{R}^2)$. Actually, by
density results we just have to prove the Brezis-Gallou\"et
inequality for smooth functions. We now chose $\phi(x)=1-|u(x)|^2$
in (\ref{14}). On the other hand, we have, setting $E=E(u)$ for the
sake of simplicity

\begin{equation}\label{11}
||\phi||^2_{L^2(\mathbb{R}^2)}\leq 4E,
\end{equation}

\begin{equation}\label{12}
||\nabla \phi||_{L^2(\mathbb{R}^2)}\leq 2||u\nabla
u||_{L^2(\mathbb{R}^2)}\leq 2||u||_{L^\infty(\mathbb{R}^2)}\sqrt{E},
\end{equation}

\begin{equation}\label{13}
\begin{split}
||\Delta \phi||_{L^2(\mathbb{R}^2)}\leq 2 ||\nabla
u||^2_{L^4(\mathbb{R}^2)}+2||u||_{L^\infty(\mathbb{R}^2)}||\Delta
u||_{L^2(\mathbb{R}^2)}\\
\leq c\left( ||\nabla u||_{L^2(\mathbb{R}^2)}||\Delta
u||_{L^2(\mathbb{R}^2)}+||u||_{L^\infty(\mathbb{R}^2)}||\Delta
u||_{L^2(\mathbb{R}^2)} \right)\\
\leq c ||\Delta u||_{L^2(\mathbb{R}^2)}(\sqrt{E}+
||u||_{L^\infty(\mathbb{R}^2)} ),
\end{split}
\end{equation}
\noindent due to classical Gagliardo-Niremberg inequality.

\noindent We then infer from (\ref{14})-(\ref{13}) that

\begin{equation}\label{16}
\begin{split}
-1+ ||u||^2_{L^\infty(\mathbb{R}^2)}\leq
c\left(\sqrt{E}(1+||u||_{L^\infty(\mathbb{R}^2)})\right)
(\log(1+R^2)^{1/2}+
\\ c(||\Delta u||_{L^2}(||u||_{L^\infty(\mathbb{R}^2)}+\sqrt{E})R^{-1})
\end{split}
\end{equation}

\noindent Then either $||u||_{L^\infty(\mathbb{R}^2)}\leq
c(1+\sqrt{E}$) and Proposition 3.1 is valid. Or,
the reverse inequality $\sqrt{E}\leq c||u||_{L^\infty(\mathbb{R}^2)}$
together with (\ref{16}) implies

\begin{equation}\label{32bis}
\begin{split}
(||u||_{L^\infty(\mathbb{R}^2)}+1)(||u||_{L^\infty(\mathbb{R}^2)}-1)
\leq c\left(\sqrt{E}(1+||u||_{L^\infty(\mathbb{R}^2)})\right)(\log(1+R^2)^{1/2}+\\c(||\Delta u||_{L^2}(||u||_{L^\infty(\mathbb{R}^2)}+1)R^{-1}),
\end{split}
\end{equation} 

\noindent and then, dividing by $||u||_{L^\infty(\mathbb{R}^2)}+1$,

\begin{equation}\label{33}
||u||_{L^\infty(\mathbb{R}^2)}\leq 1+c\left(\sqrt{E}\log(1+R^2)^{1/2}
+||\Delta u||_{L^2(\mathbb{R}^2)}R^{-1} \right).
\end{equation}

\noindent Choosing $R=||\Delta u||_{L^2}+1$ completes the proof of
 the Proposition.

\subsection{The proof of Proposition \ref{global}}

Consider a solution of (\ref{gp}) that starts from $u_0$ in
$X^2(\mathbb{R}^2)$. This solution is obtained by a fixed point
argument on the Duhamel's form of the equation. Therefore it is
classical to observe that either this solution is defined for all
time, or that the $X^2$ norm of the solution blows up in finite
time. We prove below that $||\nabla u||_{H^1(\mathbb{R}^2)}$ remains
bounded. Since the $L^2$-norm of $\nabla u$ is bounded by the
energy, it is standard to observe that we just have to bound
$||\Delta u||_{L^2(\mathbb{R}^2)}$ along the trajectories. We first
seek for an upper bound to $||u_t||_{L^2(\mathbb{R}^2)}$

Set $w(t)=u_t(t)$. Then $w(t)$ is solution in $L^2(\mathbb{R}^2)$

\begin{equation}\label{19}
iw_t+ \Delta w+w(1-|u|^2)-2{\rm Re}(w\bar{u})u=0,
\end{equation}

\noindent supplemented with initial data $w(0)=i\Delta
u_0+iu_0(1-|u_0|^2)$ in $L^2(\mathbb{R}^2)$ (once again we use that
$u_0$ has finite energy). Let us multiply (\ref{19}) by $\bar{w}$
and let us integrate the imaginary part of the resulting equation.
We obtain

\begin{equation}\label{20}
\frac{1}{2}\frac{d}{dt}||w||^2_{L^2(\mathbb{R}^2)}=2\int_{\mathbb{R}^2}{\rm
Re}(w\bar{u}){\rm Im}(w\bar{u})dx \leq
2||u||^2_{L^\infty(\mathbb{R}^2)}||w||^2_{L^2(\mathbb{R}^2)}.
\end{equation}

\noindent Using now Proposition \ref{bg}, we infer from this
inequality that

\begin{equation}\label{21}
\frac{1}{2}\frac{d}{dt}||w||^2_{L^2(\mathbb{R}^2)}\leq
c\left(1+E\right)\left(1+\log(1+||\Delta
u||^2_{L^2(\mathbb{R}^2)})\right)||w||^2_{L^2(\mathbb{R}^2)}.
\end{equation}

\noindent On the other hand, going back once more to the equation
(\ref{gp}), we also have

\begin{equation}\label{22}
||\Delta u||_{L^2(\mathbb{R}^2)}\leq ||w||_{L^2(\mathbb{R}^2)} +
2||u||_{L^\infty(\mathbb{R}^2)}\sqrt{E} \leq
||w||_{L^2(\mathbb{R}^2)}+ c(1+E)\log(1+||\Delta
u||^2_{L^2(\mathbb{R}^2)})^{1/2}.
\end{equation}

\noindent This last inequality implies that there exists a constant
$C_E$ that depends on $E$ such that $||\Delta
u||_{L^2(\mathbb{R}^2)}\leq C_E(1+||w||_{L^2(\mathbb{R}^2)})$. If we
substitute this in (\ref{21}), we infer that

\begin{equation}\label{23}
\frac{1}{2}\frac{d}{dt}||w||^2_{L^2(\mathbb{R}^2)}\leq \tilde{C}_E
(1+\log(2+||w||^2_{L^2(\mathbb{R}^2)}))||w||^2_{L^2(\mathbb{R}^2)}.
\end{equation}

\noindent Therefore the exist $a,b$ depending on
$||u_0||_{X^2(\mathbb{R}^2)}$ and on $E$ such that
$||w(t)||_{L^2(\mathbb{R}^2)}\leq e^{ae^{bt}}$. At this stage, we
infer from this and from (\ref{22}) that the $L^2$-norm of $\Delta
u$ cannot blow up in finite time. Then, due to Proposition \ref{bg},
the $L^\infty$ norm of $u$ is also controlled. This completes the
proof of Proposition \ref{global}.

\vskip 1cm

 \centerline{ACKNOWLEDGEMENT} The author would like to
thank Alberto Farina and Cl\'ement Gallo for valuable discussions
about this work. Thanks also to the anonymous referees 
for helpful remarks and comments.

\vfill\break

\end{document}